\documentclass[11pt, reqno]{amsart}
\usepackage{indentfirst, amssymb, amsmath, amsthm, mathrsfs, setspace, indentfirst, enumerate,  mathrsfs, amsmath, amsthm}
\usepackage[bookmarksnumbered, colorlinks, plainpages]{hyperref}
\usepackage{mathrsfs}

\newtheorem*{theo2A}{Theorem 2.A}
\newtheorem*{theo3A}{Theorem 3.A}

\newtheorem*{Conj3A}{Conjecture 3.A}
\newtheorem{ques}{Question}[section]
\newtheorem{theo}{Theorem}[section]
\newtheorem{lem}{Lemma}[section]
\newtheorem{cor}{Corollary}[section]

\newtheorem{exm}{Example}[section]

\newtheorem{rem}{Remark}[section]

\newcommand{\be}{\begin{equation}}
\newcommand{\ee}{\end{equation}}
\newcommand{\beas}{\begin{eqnarray*}}
\newcommand{\eeas}{\end{eqnarray*}}
\newcommand{\bea}{\begin{eqnarray}}
\newcommand{\eea}{\end{eqnarray}}

\numberwithin{equation}{section}
\begin{document}
\title [E\MakeLowercase{ntire solutions of a certain type differential-difference}.....]{\LARGE
E\Large\MakeLowercase {ntire solutions of a certain type differential-difference equation and differential-difference analogue of} B\MakeLowercase{r\"{u}ck conjecture}}
\date{}
\author[J. F. X\MakeLowercase{u}, S. M\MakeLowercase{ajumder} \MakeLowercase{and} D. P\MakeLowercase{ramanik}]{J\MakeLowercase{unfeng} X\MakeLowercase{u}, S\MakeLowercase{ujoy} M\MakeLowercase{ajumder}$^*$ \MakeLowercase{and} D\MakeLowercase{ebabrata} P\MakeLowercase{ramanik}}
\address{Department of Mathematics, Wuyi University, Jiangmen 529020, Guangdong, People's Republic of China.}
\email{xujunf@gmail.com}
\address{Department of Mathematics, Raiganj University, Raiganj, West Bengal-733134, India.}
\email{sm05math@gmail.com, sjm@raiganjuniversity.ac.in}
\address{Department of Mathematics, Raiganj University, Raiganj, West Bengal-733134, India.}
\email{debumath07@gmail.com}

\renewcommand{\thefootnote}{}
\footnote{*\emph{Corresponding Author}: Sujoy Majumder.}
\renewcommand{\thefootnote}{\arabic{footnote}}
\setcounter{footnote}{0}

\begin{abstract} In the paper, we find out the precise form of the finite order entire solutions of the following differential-difference equation 
\[f^{(k)}(z)=\sideset{}{^n_{j=0}}{\sum} a_j f(z+jc),\]
where $a_0, a_1,\ldots,a_n(\neq 0)\in\mathbb{C}$. Also in the paper we study the differential-difference analogue of Br\"{u}ck conjecture and derive a uniqueness result of finite order entire function $f(z)$ having a Borel exceptional small function of $f(z)$, when $f^{(k)}(z)$ and $\sideset{}{^n_{j=0}}{\sum} a_j f(z+jc)$ share a small function of $f(z)$. The obtained results, significantly generalize and improve the results due to Liu and Dong (Some results related to complex differential-difference equations of certain types, Bull. Korean Math. Soc., 51 (5) (2014), 1453-1467). Some examples are given to ensure the necessity of the condition (s) of our main results.
\end{abstract}
\thanks{Typeset by \AmS -\LaTeX}
\maketitle

\section{{\bf Introduction}}
In the paper, we assume that the reader is familiar with standard notation and main results of Nevanlinna Theory (see \cite{WKH1, YY1}). We use notations $\rho(f)$, $\rho_1(f)$, $\mu(f)$ and $\lambda(f)$ for the order, the hyper-order, the lower order and the exponent of convergence of zeros of a meromorphic function $f$ respectively. If $\mu(f)=\rho(f)$, we say that $f$ is of regular growth.  As usual, the abbreviation CM means ``counting multiplicities'', while IM means ``ignoring multiplicities''.

\smallskip
The Borel exceptional small function $a(z)$ of $f(z)$ is defined by
\[\lambda(f-a)=\limsup\limits_{r\to\infty} \frac{\log^+ N(r,0;f-a)}{\log r}<\rho(f)\]
where $\lambda(f-a)$ is the exponent of convergence of zeros of $f(z)-a(z)$.

\smallskip
In the paper, the forward difference $\Delta^n_c f$ for each $n\in\mathbb{N}$ is defined in the standard way by $\Delta^1_c f(z)=\Delta_c f(z)=f(z+c)-f(z)$ and $\Delta_c^n f(z)=\Delta_c\left(\Delta_c^{n-1}f(z)\right)$ $n\geq 2$. Moreover $\Delta_c^n f(z)=\sum_{j=0}^n(-1)^{n-j}C^{j}_{n}f(z+jc)$, where $C^{j}_{n}$ is a combinatorial number. 

\smallskip
We introduce the linear difference operator $L_c(f)$ by
\bea\label{e2.0}L_c(f(z))=\sideset{}{^n_{j=0}}{\sum} a_j(z) f(z+jc) (\not\equiv 0),\eea
where $a_0, a_1,\ldots,a_n(\not\equiv 0)$ are entire functions. If $a_j=(-1)^{n-j}C^{j}_{n}$, then $L_c(f)=\Delta_c^n f$. If an equation includes derivatives, shifts or differences of $f$, then the equation is called differential-difference equation.

\smallskip
Meromorphic solutions of complex difference and differential-difference equations and the uniqueness problems of meromorphic functions sharing value (s) with their derivatives, shifts, different types of complex difference and differential-difference operators have become an area of current interest and the study is based on the Nevanlinna value distribution of difference operators established by Halburd and Korhonen \cite{HK1} and Chiang and Feng \cite{CF1} respectively. 

\smallskip
The paper is organized as follows. In Section 2, we will find out the precise form of the finite order entire solutions of the differential-difference equation  $f^{(k)}(z)=L_c(f(z))$, where $a_0, a_1,\ldots,a_n(\neq 0)\in\mathbb{C}$. In Section 3, we will consider the differential-difference analogue of Br\"{u}ck conjecture and derive a uniqueness result of finite order entire function $f(z)$ having a Borel exceptional small function of $f(z)$, when $f^{(k)}(z)$ and $L_c(f(z))$ share a small function of $f(z)$.

\section{\bf{The entire solution of the equation $f^{(k)}=L_c(f)$}}
The differential-difference (or delay-differential) equation $f^{(1)}(x)=f(x-k)$, where $k>0$ is well known and studied in real analysis, some results can
be found in the book \cite{BC1}. It seems that we have little understanding on entire solutions of complex differential-difference equation
\bea\label{a} f^{(1)}(z)=f(z+c),\eea
where $c\in\mathbb{C}\backslash \{0\}$. The solutions of (\ref{a}) exist, for example, $f(z)=e^z$ is a solution of (\ref{a}), where $c=2k\pi i$, $k\in\mathbb{Z}$ and $f(z)=\sin z\;\text{or}\; \cos z$ are also solutions of (\ref{a}) for suitable $c$. Clearly, the equation (\ref{a}) has no rational solutions.

\smallskip
In 2014, Liu and Dong \cite{LD1} proved that: if $f$ is a meromorphic solution of (\ref{a}), then $\rho(f)\geq 1$. In the same paper, for the entire solutions of (\ref{a}), Liu and Dong \cite{LD1} obtained the following result. 
 
\begin{theo2A}\cite[Theorem 3.5]{LD1} If $f$ is an entire solution of (\ref{a}) and $\lambda(f)<\sigma(f)<+\infty$, then $\rho(f)=1$. Furthermore, one of the following three cases holds: 
\begin{enumerate}
\item[(i)] $f(z)=(b_1z+b_0)e^{ez+B}$, where $b_1\in\mathbb{C}\backslash \{0\}$ and $c=\frac{1}{e}$;
\item[(ii)] $f(z)=b_0e^{Az+B}$, where $c=\frac{ln|A|+i(\arg A+2ki\pi )}{A}$ and $A\in\mathbb{C}\backslash \{0\}$;
\item[(iii)] $f(z)=g(z)e^{Az+B}$, where $A\in\mathbb{C}\backslash \{0\}$ and $g(z)$ is a transcendental entire function such that $g^{(1)}(z)=A[g(z+c)-g(z)]$ and $\rho(g)<1$.
\end{enumerate}
\end{theo2A}

\smallskip
In the section, we consider the following differential-difference equation:
\bea\label{xx} f^{(k)}=L_c(f),\eea
where $a_0,a_1,\ldots,a_n(\neq 0)\in\mathbb{C}$.
The solution of Eq. (\ref{xx}) exist, for example, $f(z)=e^z$ is a solution of Eq. (\ref{xx}), where
$c=2\pi i$ and $a_0+a_1+\ldots+a_n=1$. Obviously, Eq. (\ref{xx}) has no rational solutions.

\smallskip
Now, we want to improve Theorem 2.A in the following directions:
\begin{enumerate}
\item[(i)] We replace $f^{(1)}(z)$ by general derivative $f^{(k)}(z)$.
\item[(ii)] We replace $f(z+c)$ by more general linear difference operator $L_c(f(z))$.
\end{enumerate}

We state our result.
\begin{theo}\label{t2.1} Let $f$ be an entire solution of (\ref{xx}) such that $\lambda(f)<\rho(f)<+\infty$ and let $a_0, a_1,\ldots,a_n(\neq 0)\in\mathbb{C}$. Then $\rho(f)=1$ and one of the following cases holds: 
\begin{enumerate}
\item[(1)] $f(z)=dH(z)e^{c_1z}$, where $c_1,d\in\mathbb{C}\backslash \{0\}$ and $H$ is a non-constant polynomial such that
\bea\label{XMP}\sideset{}{_{j=0}^n}{\sum} a_je^{jcc_1}H(z+jc)=\sideset{}{_{j=0}^k}{\sum} C^k_j c_1^{k-j}H^{(j)}(z) \;\text{and}\;c_1^k=\sideset{}{_{j=0}^n}{\sum} a_je^{jcc_1}\in\mathbb{C}\backslash \{0\};\eea
\item[(2)] $f(z)=dH(z)e^{c_1z}$, where $d\in\mathbb{C}\backslash \{0\}$ and $H$ is a transcendental entire function such that $\rho(H)<1$ and (\ref{XMP}) holds.
\end{enumerate}
\end{theo}

\smallskip
In order to prove Theorem \ref{t2.1}, we need the following lemmas.

\begin{lem}\label{l3}\cite[Theorem 1.51]{YY1} Suppose that $f_1, f_2,\ldots, f_n\; (n\geq 2)$ are meromorphic functions and $g_1, g_2, \ldots, g_n$ are entire functions satisfying the following conditions 
\begin{enumerate}
\item[(i)] $\sum_{j=1}^nf_je^{g_j}=0$;
\item[(ii)] $g_i-g_j$ are non-constants for $1\leq i<j\leq n$;
\item[(iii)] $T(r,f_j)=o\left(T(r,e^{g_h-g_k})\right)$ $(r\rightarrow \infty,r\not\in E)$ for $1\leq j\leq n$, $1\leq h < k\leq n$.
\end{enumerate}
Then $f_j\equiv 0$ for $j=1,2,\ldots,n$.
\end{lem}

\begin{lem}\label{l2} \cite[Corollary 2]{GGG1} Let $f$ be transcendental meromorphic function of finite order $\rho(f)$, let $\Gamma=\{(k_1,j_1),(k_2,j_2),\ldots,(k_m,j_m)\}$ denote a finite set of distinct pairs of integers that satisfy $k_i>j_i\geq 0$ for $i=1,2,\ldots,m$ and let $\varepsilon>0$ be a given constant. Then, there exists a subset $E\subset (1,\infty)$ that has a finite logarithmic measure, such that for all $z$ satisfying $|z|\not\in E\cup [0,1]$ and for all $(k,j)\in\Gamma$, we have
\beas \left|f^{(k)}(z)/f^{(j)}(z)\right|\leq |z|^{(k-j)(\rho(f)-1+\varepsilon)}.\eeas
\end{lem}

\begin{lem}\label{em1}\cite[Lemma 3.3]{BL1} Let $f$ be a transcendental meromorphic function such that $\rho(f)<1$. Let $h>0$. Then there exists an $\varepsilon$-set $E$ such that
\[f^{(1)}(z+c)/f(z+c)\rightarrow 0\;\;\text{and}\;\;f(z+c)/f(z)\rightarrow 1\;\text{as}\; z\rightarrow \infty\;in\;\mathbb{C}\backslash  E\] 
uniformly in $c$ for $|c|\leq h$. Further $E$ may be chosen so that for large $z$ not in $E$ the function $f$ has no zeros or poles in $|\zeta -z|\leq h$.
\end{lem}

\begin{lem}\label{l4} \cite[Corollary 2.6]{CF1} Let $f$ be a meromorphic function of finite order $\rho$ and let $\eta_1,\;\eta_2\in\mathbb{C}$ such that $\eta_1\neq \eta_2$. Then for any $\varepsilon>0$, we have
\beas m\left(r,f(z+\eta_1)/f(z+\eta_2)\right)=O(r^{\rho-1+\varepsilon}).\eeas
\end{lem}

\begin{lem}\label{l6}\cite[Theorem 2.3]{LY1} Let $f$ be a transcendental meromorphic solution of finite order $\rho$ of a difference equation of the form
\[U(z,f)P(z,f)=Q(z,f),\]
where $U(z,f),\; P(z,f)$ and $Q(z,f)$ are difference polynomials such that the total degree $\deg\left(U(z,f)\right)=n$ in $f$ and its shifts and $\deg\left(Q(z,f)\right)\leq n$. Moreover, we assume that $U(z,f)$ contains just one term of maximal total degree in $f$ and its shifts. Then for each $\varepsilon>0$, 
\[m(r,P(z,f))=O(r^{\rho-1+\varepsilon })+S(r,f)\]
possible outside of an exceptional set of finite logarithmic measure.
\end{lem}

\begin{rem}\label{r1} From the proof of Lemma \ref{l6}, we see that if the coefficients of
$U(z,f), P(z,f),\\ Q(z,f)$, say $a_{\lambda}$ satisfy $m(r,a_{\lambda})=S(r,f)$, then the same conclusion still holds.
\end{rem}

\begin{lem}\label{l2.1} Let $f$ be a finite order transcendental entire function such that $\rho(f)>1$ and $\lambda(f-a)<\rho(f)$, where $a$ is an entire function with $\rho(a)<\rho(f)$. If $f$ is a solution of the equation 
\bea\label{e2.1}\label{e2.2} f^{(k)}-b=\left(L_c(f)-d\right)e^{Q},\eea
where $k, n\in\mathbb{N}$, $c\in\mathbb{C}$, $Q$ is a polynomial and $a_j$'s are entire functions such that $\rho(a_j)<\rho(f)-1\;(j=0,\ldots, n)$, $a_n\not\equiv 0$ and $b,\;d$ be small function of $f$, then $\deg(Q)=\rho(f)-1$.
\end{lem}

\begin{rem} The conditions ``$\rho(f)>1$'' and ``$\rho(a_j)<\rho(f)-1$'' in Lemma \ref{l2.1} are necessary as shown in the following examples.
\end{rem}

\begin{exm} Let $a_1(z)=ze^{-2z+1},\;a_2(z)=ze^{-4z-2}$ and $Q(z)=-2$. Clearly $f(z)=e^{z^2}$ is a solution of the equation $f^{(1)}(z)=\left(a_1(z)f(z+1)+a_2(z)f(z+2)\right)e^{Q(z)}$. Note that $\rho(a_{i})=1=\rho(f)-1$ for $i=1,2$, but $\deg(Q)\neq \rho(f)-1$.
\end{exm}
\begin{exm}Let $a_1(z)=z^2e^{-3z^2+1},\;a_2(z)=2z^2e^{-3z^2-9z-6}$ and $Q(z)=-3z-2$. Clearly $f(z)=e^{z^3}$ is a solution of the equation 
$f^{(1)}(z)=(a_1(z)f(z+1)+a_2(z)f(z+2))e^{Q(z)}$. Note that $\rho(a_{j})=2=\rho(f)-1$ for $i=1,2$, but $\deg(Q)\neq \rho(f)-1$.
\end{exm}

\smallskip
\begin{proof}[{\bf Proof of Lemma \ref{l2.1}}]
By the given conditions, we have $\rho(a)<\rho(f)$ and $\lambda(f-a)<\rho(f)$. Then there exist an entire function $H(\not\equiv 0)$ and a polynomial $P$ such that 
\bea\label{e.1} f=a+He^{P},\eea
where $\lambda(H)=\rho(H)<\rho(f-a)$ and $\deg(P)=\rho(f-a)$. Since $\rho(a)<\rho(f)$, we have $\rho(f-a)=\rho(f)$. We know that $\mu(f)=\rho(f)=\deg(P)>1$. Since $\rho(a_j)<\rho(f)-1<\rho(f)=\mu(f)$, by Theorem 1.18 \cite{YY1}, we get $T(r,a_j)=S(r,f)$ for $j=0,1,\ldots,n$. Similarly $T(r,H)=S(r,f)$. Now from (\ref{e.1}), we have
\bea\label{e.2}\label{al.2a} L_c(f(z))=\sideset{}{^n_{j=0}}{\sum}a_j(z)a(z+jc)+H_n(z)e^{P(z)},\eea
where 
\bea\label{e.2a}\label{a1.a1} H_n(z)=\sideset{}{^n_{j=0}}{\sum} a_j(z)H(z+jc)e^{P(z+jc)-P(z)}.\eea

Also from (\ref{e.1}), we get 
\bea\label{e.e1} f^{(k)}= a^{(k)}+ P_k(H)e^{P},\eea
where 
\bea\label{e.2aa} P_k(H)=H^{(k)}+p_{k-1}H^{(k-1)}+\ldots + p_0H\eea 
and $p_i$'s are polynomials such that 
\bea\label{e.2aaa}\deg(p_i)=(k-i)\deg(P^{(1)})\eea
for $i=0,1,\ldots, k-1$. Now from (\ref{e2.2})-(\ref{e.2}) and (\ref{e.e1}), we have 
\bea\label{e.3} a^{(k)}-b+P_k(H)e^{P}=\Big(\sideset{}{^n_{j=0}}{\sum}a_j(z)a(z+jc)-d+H_ne^{P}\Big)e^{Q}.\eea

Let $w_1=P_k(H)$, $w_2=a^{(k)}-b$, $w_3=H_n$ and $w_4=\sum_{j=0}^na_j(z)a(z+jc)-d$. It is clear that $T(r,w_j)=S(r,f)$ for $j=1,2,3,4$. Then from (\ref{e.3}), we have 
\bea\label{e.4} w_1e^P+w_2=(w_3e^P+w_4)e^Q.\eea
 
Obviously $\deg(Q)\leq \deg(P)$. Let $P(z)=\tilde{a}_sz^s+\tilde{a}_{s-1}z^{s-1}+\ldots+\tilde{a}_0$. Clearly $s\geq 2$.

Now we divide following two cases.\par

\medskip
{\bf Case 1.} Suppose $\deg(Q)=\deg(P)$. Let $Q(z)=\tilde{b}_sz^s+\tilde{b}_{s-1}z^{s-1}+\cdots+\tilde{b}_0\;(\tilde{b}_s\neq 0)$. There are three possibilities: $(1)\;\tilde{a}_s=\tilde{b}_s$, $(2)\;\tilde{a}_s=-\tilde{b}_s$ and $(3)\;\tilde{a}_s\neq \tilde{b}_s\;\text{and}\;\tilde{a}_s\neq -\tilde{b}_s$. Therefore we consider following three sub-cases.\par

\medskip
{\bf Sub-case 1.1.} Suppose $\tilde{a}_s=\tilde{b}_s$. Then from (\ref{e.4}), we have
\bea\label{e.5} (w_1-w_4e^{Q-P})e^P-w_3e^{P+Q}+w_2=0.\eea

Clearly $T(r,e^{Q-P})=o(T(r,e^P))$ and $T(r,e^{Q-P})=o(T(r,e^Q))$. Also $T(r,w_3)=o(T(r,e^P))$ and $T(r,w_3)=o(T(r,e^Q))$. Consequently
\[T(r,w_1-w_4e^{Q-P})=o(T(r,e^P)), T(r,w_1-w_4e^{Q-P})=o(T(r,e^Q)),\ldots, T(r,w_2)=o(T(r,e^{P+Q})).\]

Therefore using Lemma \ref{l3} to (\ref{e.5}), we get $w_1-w_4e^{Q-P}\equiv 0$, $w_2\equiv 0$ and $w_3\equiv 0$.  Now $w_3\equiv 0$ implies that
\bea\label{e.6} \sideset{}{^n_{j=1}}{\sum}e^{P(z+jc)-P(z)}a_j(z)H(z+jc)/H(z)=-a_0(z).\eea

Let $R(z)=P(z+c)-P(z)$. By the given condition, we have $\rho(a_j)<\rho(f)-1$, $j=0,1,\ldots,n$ and so $\rho(a_j)<\rho(e^R)=\mu(e^R)$. Clearly $T(r,a_j)=S(r,e^R)$, $j=0,1,\ldots,n$. Also by the definition of order, there exists a $R>0$ such that $T(r,a_j)=O(r^{\rho(a_j)+\varepsilon})$ holds for $r\geq R$, where $\varepsilon >0$ is arbitrary.  
Now from Lemma \ref{l4}, we have 
\bea\label{e.6a} m\left(r,H(z+jc)/H(z)\right)+m\left(r,H(z)/H(z+jc)\right)=O\left(r^{\rho (H)-1+\varepsilon}\right).\eea

We consider following sub-cases.

\medskip
{\bf Sub-case 1.1.1.} Let $n=1$. Then from (\ref{e.6}), we have
\bea\label{ee.1} e^{R(z)}a_1(z)H(z+c)/H(z)=-a_0(z).\eea

By the given condition $a_1\not=0$. If $a_0\equiv 0$, then from (\ref{ee.1}) we get $H\equiv 0$, which is impossible. Hence $a_0\not\equiv 0$. 
Now (\ref{ee.1}) shows that $a_1(z)H(z+c)/H(z)$ is an entire function. Note that
\[T\left(r,a_1(z)H(z+c)/H(z)\right)=m\left(r,a_1(z)H(z+c)/H(z)\right)\leq m(r,a_1)+m\left(r,H(z+c)/H(z)\right)\]
and so
\bea\label{ee.2}T\left(r,a_1(z)H(z+c)/H(z)\right)=m\left(r,a_1(z)H(z+c)/H(z)\right)=O\left(r^{\sigma+\varepsilon}\right),\eea
where $\sigma=\max\{\rho(H)-1, \rho(a_1)\}$. 

Since $\rho(H)<\rho(f)$ and $\rho(a_1)<\rho(f)-1$, we have $\sigma<\rho(f)-1=\rho(e^R)=\mu(e^R)$. Now (\ref{ee.2}) yields
$\rho\big(a_1(z)H(z+c)/H(z)\big)=\sigma<\mu(e^{R})$ and so
\[T\big(r,a_1(z)H(z+c)/H(z)\big)=S(r,e^{R}).\]

Then from (\ref{ee.1}), we get a contradiction.

\medskip
{\bf Sub-case 1.1.2.} Let $n\geq 2$.

First we suppose that $a_0\not\equiv 0$. 

Clearly $Q_1(z+(j-1)c)\ldots Q_1(z+c)=e^{P(z+jc)-P(z+c)}$, $j=2,3,\ldots,n$, where $Q_1=e^{R}$. Then (\ref{e.6}) can be written as
\bea\label{e.7} U(z,Q_1)Q_1=-a_0,\eea
where 
\beas U(z,Q_1(z))=\sum\limits_{1\leq j\leq n}a_{j}(z)Q_1(z+(j-1)c)Q_1(z+(j-2)c)\cdots Q_1(z+c)H(z+jc)/H(z).\eeas

Now from (\ref{e.7}), we see $\deg(U(z,Q_1))=n-1\geq 1$. 
Also from above, we have
\bea\label{ee.3} m(r,a_j(z)H(z+jc)/H(z))=O(r^{\sigma+\varepsilon}),\eea
where $\sigma=\max\{\rho(H)-1, \rho(a_j)\}$, $j=0,1,\ldots,n$. 

Since $\rho(H)<\rho(f)$ and $\rho(a_j)<\rho(f)-1$, it follows that $\sigma<\rho(f)-1=\rho(e^R)=\mu(e^R)$. Let us 
choose $\varepsilon>0$ such that $\sigma+2\varepsilon <\rho(e^R)=\mu(e^R)$. By the definition of lower order, there exists $R>0$ such that 
\[T(r,e^R)>r^{\mu(e^R)-\varepsilon}=r^{\rho(e^R)-\varepsilon}\]
holds for $r\geq R$. Then from (\ref{ee.3}), we get
\[\lim\limits_{r\to\infty}\frac{m(r,a_j(z)H(z+jc)/H(z))}{T(r,e^R)}=0\]
and so 
\bea\label{ee.4} m\left(r,a_j(z)H(z+jc)/H(z)\right)=S(r,e^{R(z)})=S(r,Q_1(z))\eea
for $j=0,1,\ldots,n$. This shows that if $a_{\lambda}$ is a coefficient of $U(z,Q_1)$, then $m(r,a_\lambda)= S(r,Q_1)$. Then in view of Remark \ref{r1} and using Lemma \ref{l6} to (\ref{e.7}), we get 
$m(r,Q_1)=S(r,Q_1)$ and so $T(r,Q_1)=S(r,Q_1)$, which is a contradiction.

\medskip
Next we suppose $a_0=0$. For the sake of simplicity we may assume that $a_1\not\equiv 0$. 
So from (\ref{e.6}) we get
\beas\label{ea.1}\sideset{}{^n_{j=2}}{\sum}a_j(z)e^{P(z+jc)-P(z+c)}H(z+jc)/H(z)=-a_1(z).\eeas

Now proceeding in the same way as done above, we get a contradiction.\par

\medskip
{\bf Sub-case 1.2.} Suppose $\tilde{a}_s=-\tilde{b}_s$. Then from (\ref{e.4}), we have
\bea\label{e.8} w_1e^P-w_4e^Q-w_3e^{P+Q}+w_2=0.\eea

Clearly $T(r,w_1)=o(T(r,e^Q))$, $T(r,w_1)=o(T(r,e^P)),\ldots$, $T(r,w_2-w_3e^{P+Q})=o(T(r,e^Q)).$
Then using Lemma \ref{l6} to (\ref{e.8}), we get $w_1\equiv 0$, $w_2-w_3e^{P+Q}\equiv 0$ and $w_4\equiv 0$.
Since $w_1\equiv 0$, we have $H\equiv 0$, which is impossible.\par

\medskip
{\bf Sub-case 1.3.} Suppose $\tilde{a}_s\neq \tilde{b}_s$ and $\tilde{a}_s\neq -\tilde{b}_s$. Then from (\ref{e.4}), we have
\bea\label{e.9} w_1e^P-w_4e^Q-w_3e^{P+Q}+w_2=0.\eea

Obviously $T(r,w_1)=o(T(r,e^P))$, $T(r,w_1)=o(T(r,e^Q)), \ldots$, $T(r,w_4)=o(T(r,e^{P+Q}).$
Then using Lemma \ref{l6} to (\ref{e.8}), we get $w_1\equiv 0$, $w_2-w_3e^{P+Q}\equiv 0$ and $w_4\equiv 0$.
Since $w_1\equiv 0$, we have $H\equiv 0$, which is impossible.\par

\medskip
{\bf Case 2.} Suppose $\deg(Q)<\deg(P)$. Then from (\ref{e.4}), we have 
\bea\label{e.10} (w_1-w_3e^Q)e^P+(w_2-w_4e^Q)=0.\eea

Then using Lemma \ref{l6} to (\ref{e.10}), we get
$w_1-w_3e^Q\equiv 0$ and $w_2-w_4e^Q\equiv 0$. Therefore $w_1-w_3e^Q\equiv 0$ implies that
\bea\label{e.11}\sideset{}{^n_{j=0}}{\sum}a_j(z)e^{P(z+jc)-P(z)}H(z+jc)/H(z)=e^{-Q(z)}P_k(H(z)/H(z).\eea

Now for $0\leq j\leq n$, we have 
$e^{P(z+jc)-P(z)}=e^{jsca_sz^{s-1}}e^{P_j(z)}=h^j(z)e^{P_j(z)}$,
where $h(z)=e^{sca_sz^{s-1}}$ and $\deg(P_j)\leq s-2$. 
Note that 
\[\rho(e^{P_j})=\deg(P_j)\leq s-2<s-1=\rho(h)=\mu(h)\]
and so $T(r, e^{P_j})=S(r,h)$ and so $T(r,e^{P_j})=S(r,h^j)$ for $j=0,1,2,\ldots, n$.
Also from above, one can easily prove that 
\bea\label{al.9} m\left(r,a_j(z)H(z+jc)/H(z)\right)=S(r,h)\eea
for $j=0,1,\ldots,n$.
Similarly we can prove that
\bea\label{al.10} m\left(r,a_j(z)H(z)/H(z+jc)\right)=S(r,h)\eea
for $j=0,1,\ldots,n$. Let
\bea\label{al.11} b_{n-j}(z)=e^{P_j(z)}a_j(z)H(z+jc)/H(z),\eea
for $j=0,1,2,\ldots, n$ and  
\bea\label{al.12} F_n(h)=\sideset{}{^n_{j=0}}{\sum}b_{n-j}h^j.\eea

Since $T(r, e^{P_j})=S(r,h)$, $j=0,1,\ldots,n$, from (\ref{al.9})-(\ref{al.11}), we get 
\[m(r,b_j)+m(r,1/b_j)=S(r,h),\; j=0,1,\ldots,n.\]

Now from the proof of Theorem 2.1 \cite{LLLX1}, we conclude that
\bea\label{e.14} n\;m(r,h)+S(r,h)=m(r,F_n(h))+S(r,h).\eea

Also from (\ref{e.11}), (\ref{al.11}) and (\ref{al.12}), we get $F_n(h)=\frac{P_k(H)}{H}e^{-Q}$ and so from (\ref{e.14}) we get
\bea\label{e.15} n\;m(r,h)+S(r,h)=m\left(r,P_k(H)e^{-Q}/H\right)\leq m\left(r,P_k(H)/H\right) + m(r,e^{-Q}).\eea

We know that $m(r,P_k(H)/H) = O(\log r)$ as $r\to\infty$. Since $h$ is a transcendental entire function, it follows that
\[\frac{m(r,P_k(H)/H)}{T(r, h)} \rightarrow 0\]
as $r\rightarrow \infty$ and so 
\bea\label{e.111} m(r,P_k(H)/H) = S(r,h).\eea

Now from (\ref{e.15}) and (\ref{e.111}), we have
\bea\label{e.112} n\;m(r,h)+S(r,h)\leq m(r,e^{-Q}).\eea

We claim that $\deg(Q)=\rho(f)-1=s-1$. If not, suppose $\deg(Q)<\rho(f)-1=s-1$. Then $\rho(e^{-Q})<s-1=\rho(h)=\mu(h)$ and so $T(r,e^{-Q})=S(r,h)$, i.e., $m(r,e^{-Q})=S(r,h)$. Now from (\ref{e.112}), we get $n\;m(r,h)=S(r,h)$, i.e., $n\;T(r,h)=S(r,h)$, which is impossible. Hence $\deg(Q)=s-1$, i.e., $\deg(Q)=\rho(f)-1$.

This completes the proof.
\end{proof}

\smallskip
\begin{proof}[{\bf Proof of Theorem \ref{t2.1}}]
Let $f$ be an entire solution of (\ref{xx}) such that $\lambda(f)<\rho(f)<+\infty$. Then there exist an entire function $H(\not\equiv 0)$ and a polynomial $P$ such that 
\bea\label{se.1} f=He^{P},\eea
where $\lambda(H)=\rho(H)<\rho(f)$ and $\deg(P)=\rho(f)$. We consider following two cases.

\medskip
{\bf Case 1.} Let $\rho(f)>1$. If we take $b=d=0$ and $Q\equiv 0$, then by Lemma \ref{l2.1}, we get $0=\deg(Q)=\rho(f)-1$, i.e., $\rho(f)=1$, which is impossible.

\medskip
{\bf Case 2.} Let $\rho(f)\leq 1$. Since $0$ and $\infty$ are Borel exceptional values of $f$, we have $\rho(f)=1$. Clearly $\rho(H)<1$ and $\deg(P)=1$. Let $P(z)=c_1z+c_0$, $c_1\neq 0$. Now from (\ref{se.1}), we have
\bea\label{se.2} L_c(f(z))=H_n(z)e^{P(z)}\quad \text{and}\quad f^{(k)}=P_k(H)e^{P},\eea
where 
\bea\label{se.3} H_n(z)=\sideset{}{^n_{j=0}}{\sum} a_jH(z+jc)e^{jcc_1}\quad\text{and}\quad P_k(H)=\sideset{}{_{j=0}^k}{\sum} C^k_j c_1^{k-j}H^{(j)}(z).\eea

Therefore using (\ref{se.2}) and (\ref{se.3}) to (\ref{xx}), we get
\bea\label{se.4} \sideset{}{^n_{j=0}}{\sum} a_je^{jcc_1}H(z+jc)=\sideset{}{^k_{j=0}}{\sum} C^k_j c_1^{k-j}H^{(j)}(z).\eea 

\medskip
First we suppose $H$ is a polynomial. Then by a simple calculation, we deduce from (\ref{se.4}) that $c_1^k=\sum_{j=0}^n a_j e^{jcc_1}$. 

\medskip
Next we suppose $H$ is transcendental. Then by Lemma \ref{l2}, for any given $\varepsilon>0\;(\rho(H)+\varepsilon<1)$, there exists a subset $E\subset (1,\infty)$ of finite logarithmic measure such that for all $z$ satisfying $|z|=r\not\in E\cup [0,1]$, we have
\bea\label{se.5} |H^{(j)}(z)/H(z)|\leq |z|^{j(\rho(H)-1+\varepsilon)},\;\;j=1,2,\ldots,k.\eea

Again by Lemma \ref{em1}, we obtain 
\bea\label{se.6}\sideset{}{_{j=0}^n}{\sum}a_je^{jcc_1}H(z+jc)/H(z)\to \sideset{}{_{j=0}^n}{\sum} a_je^{jcc_1}.\eea

Now using (\ref{se.5}) and (\ref{se.6}) to (\ref{se.4}), we get $\sum_{j=0}^n a_je^{jcc_1}=c_1^k$. Finally we have 
\[f(z)=dH(z)e^{c_1z},\]
where $d\in\mathbb{C}\backslash \{0\}$ and $H$ is a transcendental entire function such that $\rho(H)<1$,
\[\sideset{}{_{j=0}^n}{\sum} a_je^{jcc_1}H(z+jc)=\sideset{}{_{j=0}^k}{\sum} C^k_j c_1^{k-j}H^{(j)}(z) \;\text{and}\;c_1^k=\sideset{}{_{j=0}^n}{\sum} a_je^{jcc_1}\in\mathbb{C}\backslash \{0\}.\]

This completes the proof.
\end{proof}

\section {{\bf Differential-difference analogue of Br\"{u}ck conjecture}}
We recall the following conjecture proposed by Br\"{u}ck \cite{RB1}.
\begin{Conj3A} Let $f$ be a non-constant entire function such that $\rho_1(f)\not\in\mathbb{N}\cup\{\infty\}$. If $f$ and $f^{(1)}$ share one finite value $a$ CM, then $f^{(1)}-a=c(f-a),\;\text{where}\;\;c\in\mathbb{C}\backslash \{0\}$. 
\end{Conj3A}

\smallskip
Though the conjecture is not settled in its full generality, it gives rise to a long course of research on the uniqueness of entire and meromorphic functions sharing a single value with its derivatives.

Recently, many authors have started to consider the sharing values problems of meromorphic functions with their difference operators or shifts. Some results on difference analogues of Br\"{u}ck conjecture were considered in \cite{CC1}-\cite{ZXC1}, \cite{HK1}, \cite{HKLRZ}-\cite{LLLX1}, \cite{MD1}, \cite{MP1}, \cite{MSP1}, \cite{QLY1}.
 
\smallskip
Br\"{u}ck \cite{RB1} proved that if $f$ is a transcendental entire function and $f$ and $f^{(1)}$ share $0$ CM, then $f^{(1)}=Af$.  In 2009, Heittokangas et al. \cite{HKLRZ} proved that if $f(z)$ is a transcendental meromorphic function with $\rho(f)<2$ and $f(z)$ and $f(z+c)$ share $0$ CM, then $f(z)=Af(z+c)$. We now consider the differential-difference analogue of Br\"{u}ck conjecture. Note that if $f^{(1)}(z)$ and $f(z+c)$ share $0$ CM, the equation $f^{(1)}(z)=Af(z+c)$ is not valid in general. As for example if $f(z)=e^{e^z}$ and $c=2\pi i$, then $f^{(1)}(z)=e^z f(z+c)$. However, for finite order entire functions the differential-difference analogue of Br\"{u}ck conjecture may hold. For example, let $f(z)=e^z$ and $e^c=2$. Then $f^{(1)}(z)$ and $f(z +c)$ share $0$ CM and $f^{(1)}(z)=\frac{1}{2}f(z+c)$. In 2014, Liu and Dong \cite{LD1} first considered the differential-difference analogue of Br\"{u}ck conjecture for finite order transcendental entire function. Now we recall their result.

\begin{theo3A}\cite[Theorem 4.3]{LD1} Let $f$ be a finite order transcendental entire function and $a$ be a Picard exceptional value of $f$. If $f^{(1)}(z)$ and $f(z+c)$ share $b(\in\mathbb{C})$ CM, then $f^{(1)}(z)-b=A(f(z+c)-b)$, where $A\in\mathbb{C}\backslash \{0\}$. Furthermore, if $b\neq 0$, then $A=\frac{b}{b-a}$. 
\end{theo3A}

In this section, we will improve Theorem 3.A in the following directions:
\begin{enumerate}
\item[(1)] We replace $f^{(1)}(z)$ by $f^{(k)}(z)$ and $f(z+c)$ by $L_c(f(z))$.
\item[(2)] We consider $a$ as an entire function such that $\lambda(f-a)<\rho(f)$ and $\rho(a)<\max\{1,\rho(f)\}$.
\item[(3)] We consider $b$ as an entire function such that $\rho(b)<\max\{1,\rho(f)\}$.
\end{enumerate}

We now state our result.

\begin{theo}\label{t1} Let $f$ be a finite order transcendental entire function and let $a$ and $b$ be entire functions such that $\lambda(f-a)<\rho(f)$, $\rho(a)\neq \rho(f)$ and $\rho(a), \rho(b)<\max\{1,\rho(f)\}$. Let $L_c(f)$ be defined as in (\ref{e2.0}), where $a_0, a_1,\ldots,a_n(\neq 0)\in\mathbb{C}$. If $f^{(k)}$ and $L_c(f)$ share $b$ CM, then
\begin{enumerate}

\smallskip
\item[(1)] $f^{(k)}-b=A(L_c(f)-b)$, where $A\in\mathbb{C}\backslash \{0\}$, $a^{(k)}-A\sideset{}{_{j=0}^n}{\sum} a_ja(z+jc)=(1-A)b$ and one of the following cases holds:

\smallskip
\begin{enumerate}
\item[(1A)] $f(z)=a(z)+dH(z)e^{c_1z}$, where $c_1,d\in\mathbb{C}\backslash \{0\}$ and $H$ is a polynomial such that 
\bea\label{XMP1}A\sideset{}{_{j=0}^n}{\sum} a_je^{jcc_1}H(z+jc)=\sideset{}{_{j=0}^k}{\sum} C^k_j c_1^{k-j} H^{(j)}(z)\;\;\text{and}\;\;c_1^k=A\sideset{}{_{j=0}^n}{\sum} a_je^{jcc_1};\eea

\smallskip
\item[(1B)] $f(z)=a(z)+dH(z)e^{c_1z}$, where $d\in\mathbb{C}\backslash \{0\}$ and $H$ is a transcendental entire function such that $\rho(H)<1$ and (\ref{XMP1}) holds;
\end{enumerate}

\smallskip
\item[(2)] $f^{(k)}-b=e^{c_1z+d_0}(L_c(f)-b)$, where $c_1(\neq 0),d_0\in\mathbb{C}$, $a\neq b=a^{(k)}\not\equiv 0$, $a_j\neq 0$ for at least one $j=0,\cdots, n-1$ such that $\sum_{j=0}^na_je^{jcc_1}=0$ and one of the following cases holds:
\begin{enumerate}

\smallskip
\item[(2A)] $f(z)=a(z)+H(z)e^{c_1z+c_0}$, where $a(\not\equiv 0)$ and $H(\not\equiv 0)$ are polynomials such that $a^{(k)}=\sum_{j=0}^n a_j a(z+jc)-e^{c_0-d_0}\sum_{j=0}^k C^k_j c_1^{k-j}H^{(j)}$;

\smallskip
\item[(2B)] $f(z)=a(z)+H(z)e^{c_1z+c_0}$, where $a$ is a transcendental entire function such that $\rho(a)<1$ and $H(\not\equiv 0)$ is a polynomial such that $a^{(k)}=\sum_{j=0}^n a_j a(z+jc)-e^{c_0-d_0}\sum_{j=0}^k C^k_j c_1^{k-j}H^{(j)}$ and $e^{cc_1}=1$;
\end{enumerate}

\smallskip
\item[(3)] $f^{(k)}-b=e^{Q}(a_1f(z+c)-b)$ and $f(z)=a+He^{P}$, where $P$ and $Q$ are polynomials and $H$ is a transcendental entire function such that $\deg(P)\geq 2$, $\deg(Q)=\deg(P)-1$, $\rho(H)\geq 1$ and $(a^{(k)}-b)e^{-Q}=a_1a(z+c)-b$.
\end{enumerate}
\end{theo}

We now observe from Theorem \ref{t1} that if one of the following conditions holds:
\begin{enumerate}
\item[(i)] $n=1$ and $a_0\not=0$;
\item[(ii)] $a^{(k)}\not=b$ and $a^{(k)}(z)=a_1a(z+c)$, when $n=1$, $a_0=0$ and $\lambda(f-a)\geq 1$,
\end{enumerate}
then the conclusion $(3)$ does not hold.

\smallskip
If $a$ and $b$ are constants, then from Theorem \ref{t1} we get the following corollary.
\begin{cor}\label{c31} Let $f$ be a transcendental entire function of finite order and let $a$ and $b$ be two constants such that $\lambda(f-a)<\rho(f)$. Let $L_c(f)$ be defined as in (\ref{e2.0}), where $a_0, a_1,\ldots,a_n(\neq 0)\in\mathbb{C}$. If $f^{(k)}$ and $L_c(f)$ share $b$ CM, then $f^{(k)}-b=A(L_c(f)-b)$, where $A\in\mathbb{C}\backslash \{0\}$, $A\sideset{}{_{j=0}^n}{\sum} a_ja(z+jc)=(A-1)b$ and one of the following cases holds:
\begin{enumerate}
\item[(i)] $f(z)=a(z)+dH(z)e^{c_1z}$, where $c_1,d\in\mathbb{C}\backslash \{0\}$ and $H$ is a polynomial such that (\ref{XMP1}) holds;
\item[(ii)] $f(z)=a(z)+dH(z)e^{c_1z}$, where $d\in\mathbb{C}\backslash \{0\}$ and $H$ is a transcendental entire function such that $\rho(H)<1$ and (\ref{XMP1}) holds.
\end{enumerate}
\end{cor}

Obviously Corollary \ref{c31} improves as well as generalizes Theorem 3.A in a large scale.

\begin{rem} Following example shows that condition ``$\rho(f)<+\infty$'' is sharp in Theorem \ref{t1}.
\end{rem}
\begin{exm} Let $f(z)=e^{e^z}$ and let $a=0$, $b=0$ and $c=2\pi i$. Here $\lambda(f-a)=0<\rho(f)=+\infty$, $\rho(a)\neq \rho(f)$, $\rho(a), \rho(b)<\max\{1,\rho(f)\}$ and 
\[e^{z}(f(z+c)-b)=f^{(1)}(z)-b.\]

Clearly $f^{(1)}(z)$ and $f(z+c)$ share $b$ CM, but $f^{(1)}(z)-b=A(f(z+c)-b)$ is not valid, where $A\in\mathbb{C}\backslash \{0\}$.
\end{exm}

\begin{rem} Following example shows that conditions ``$\rho(a)\neq \rho(f)$'' and ``$\rho(a),\; \rho(b)<\max\{1,\rho(f)\}$'' are sharp in Theorem \ref{t1}.
\end{rem}

\begin{exm} Let $f(z)$ is a transcendental entire function such that $\rho(f)\geq 1$. Let $L_c(f(z))=f(z)+f(z+c)$, $a(z)=f(z)-z$ and $b(z)=-f(z)-f(z+c)+2f^{1)}(z)$. Note that $\lambda(f-a)=0<\rho(f)$, $\rho(a)=\rho(f)$ and $\rho(b)\leq \rho(f)$. Clearly $L_c(f)-b=2(f^{(1)}-b)$, which means that $f^{(1)}$ and $L_c(f)$ share $b$ CM. But, $f$ does not satisfy the form of $f$ in Theorem \ref{t1}.
\end{exm}

\smallskip
In the proof of Theorem \ref{t1} below, we make use of the following key lemmas. Now we recall these lemmas here.
\begin{lem}\label{l7} \cite[Corollary 8.3]{CF1} Let $\eta_1,\eta_2$ be two arbitrary complex numbers and let $f$ be a meromorphic function of finite order $\sigma$. Let $\varepsilon>0$ be given. Then there exists a subset $E\subset (0,\infty)$ with finite logarithmic measure such that for all $r\not\in E\cup[0,1]$, we have
\beas \exp(-r^{\sigma -1+\varepsilon })\leq \left|f(z+\eta_1)/f(z+\eta_2)\right|\leq \exp(r^{\sigma -1+\varepsilon }).\eeas
\end{lem}

\begin{lem}\label{l8}\cite[Corollary 3.2]{LLLX1} Let $f$ be a non-constant meromorphic solution of  
\[\sideset{}{_{i=0}^k}{\sum} b_if(z+ic)=R(z),\]
where $R$ is polynomial and $b_i\in\mathbb{C}$, $i=0,1,\ldots, k$ with $b_kb_0\neq 0$, $c\in \mathbb{C}\backslash \{0\}$, $k\in \mathbb{N}$. Then either $\rho(f)\geq 1$ or $f$ is a polynomial. If $b_k\neq \pm b_0$, then $\rho(f)\geq 1$.
\end{lem}

\begin{lem} \label{l6a} \cite[Lemma 1.3.1]{IL1} Let $P(z)=\sum_{i=1}^{n}a_{i}z^{i}\;(a_{n}\not=0)$. Then $\forall$ $\varepsilon>0$, $\exists$ $r_{0}>0$ such that $\forall$ $r=|z|>r_{0}$ the inequalities $(1-\varepsilon)|a_{n}|r^{n}\leq |P(z)|\leq (1+\varepsilon)|a_{n}|r^{n}$ hold. 
\end{lem}

\begin{lem}\label{l2a}\cite[Lemma 2]{LZY1} Let $f$ be an entire function and suppose $|f^{(k)}(z)|$ is unbounded on some ray $\arg z=\varphi$. Then there exists an infinite sequence of points $z_n=r_ne^{i\varphi}$, where $r_n\to \infty$ such that $f^{(k)}(z_n)\to \infty$ and
\[|f(z_n)|\leq (1+o(1))|z_n|^k|f^{(k)}(z_n)|\;\text{as}\;z_n\to\infty.\]
\end{lem}

\smallskip
\begin{proof}[{\bf Proof of Theorem \ref{t1}}]
By the given conditions, we get $f=a+He^{P}$, where $H(\not\equiv 0)$ is an entire function and $P$ is a polynomial such that $\lambda(H)=\rho(H)<\rho(f-a)$ and $\deg(P)=\rho(f-a)$. 

\medskip
First we suppose $\rho(f)<1$. 

Since $\rho(a)<\max\{1,\rho(f)\}$ and $\rho(a)\neq \rho(f)$, we have $\rho(a)<1$ and so $\rho(f-a)<1$. Consequently $\lambda(f-a)<\rho(f)\leq \max\{\rho(a), \rho(f)\}=\rho(f-a)$. This shows that $0$ and $\infty$ are the Borel exceptional values of $f-a$ and so $f-a$ is a function of regular growth. Therefore by Theorem 2.11 \cite{YY1}, we have $\rho(f-a)\in\mathbb{N}$, which is impossible.

\medskip
Next we suppose $\rho(f)\geq 1$. 

In this case, we have $\rho(a)<\rho(f)$ and $\rho(b)<\rho(f)$. Therefore $\rho(H)<\rho(f)$ and $\deg(P)=\rho(f)$. Let 
\bea\label{al.1a} P(z)=c_sz^s+c_{s-1}z^{s-1}+\cdots+c_0,\eea
where $c_s(\neq 0),\ldots,c_0\in\mathbb{C}$ and $s\in\mathbb{N}$. Therefore $\rho(f)=s$. As $f^{(k)}$ and $L_c(f)$ share $b$ CM, there exists a polynomial function $Q$ such that 
\bea\label{al.4a}f^{(k)}-b=\Big(\sideset{}{^n_{j=0}}{\sum}a_jf(z+jc)-b\Big)e^{Q}.\eea

Clearly (\ref{al.4a}) yields $\deg(Q)\leq \rho(f)$. Therefore using (\ref{al.2a}) and (\ref{e.e1}) to (\ref{al.4a}), we have 
\bea\label{m1} a^{(k)}-b-\Big(\sideset{}{^n_{j=0}}{\sum}a_ja(z+jc)-b\Big)e^{Q}=-\left(P_k(H)-H_ne^{Q}\right)e^{P},\eea
where $H_n$ and $P_k(H)$ are defined respectively in (\ref{a1.a1}) and (\ref{e.2aa}).

Now we divide following two cases.\par

\medskip
{\bf Case 1.} Suppose $\rho(f)<2$. Since $\deg(P)=\rho(f)$, it follows that $\deg(P)<2$ and so $\deg(P)=1$. Consequently $\rho(f)=1$. Therefore by the given conditions, we see that $\rho(H)<1$, $\rho(a)<1$, $\rho(a)\neq 1$ and $\rho(b)<1$. Since $P(z)=c_1z+c_0$, where $c_1\neq 0$, then from (\ref{e.2aa}), we get $p_i=C^k_i c_1^{k-i}$, $i=0,1,2,\ldots, k$.

Now we divide the following two sub-cases.\par

\medskip
{\bf Sub-case 1.1.} Suppose $\deg(Q)=0$. Let $e^Q=A$ and so $f^{(k)}-b=A(L_c(f)-b)$. On the other hand from (\ref{m1}), we have
\bea\label{al.6} a^{(k)}-b-A\Big(\sideset{}{^n_{j=0}}{\sum}a_ja(z+jc)-b\Big)=-\left(P_k(H)-AH_n\right)e^{P}.\eea

Note that $\rho\left(P_k(H)-AH_n\right)<\rho(f)=\rho\left(e^P\right)=\mu\left(e^P\right)$ and so $T\left(r,	P_k(H)-AH_n\right)=S\left(r,e^P\right)$. Now using Lemma \ref{l3} to (\ref{al.6}), we get $P_k(H)-AH_n=0$ and $a^{(k)}-A\sum_{j=0}^na_ja(z+jc)=(1-A)b$.
On the other hand, from (\ref{a1.a1}) and (\ref{e.2aa}), we get 
\beas\label{vv.1} A\sideset{}{^n_{j=0}}{\sum} a_je^{jcc_1}H(z+jc)=\sideset{}{^k_{j=0}}{\sum} C^k_i c_1^{k-j} H^{(j)}(z).\eeas 

If $H$ is a polynomial, then from the proof of Theorem \ref{t2.1}, we get $f(z)=a(z)+dH(z)e^{c_1z}$, where $c_1,d\in\mathbb{C}\backslash \{0\}$ such that (\ref{XMP1}) holds and $a^{(k)}-A\sideset{}{_{j=0}^n}{\sum} a_ja(z+jc)=(1-A)b$. Again if $H$ is transcendental, then proceeding in the same way as done in the proof of Theorem \ref{t2.1}, we obtain $A\sum_{j=0}^n a_je^{jcc_1}=c_1^k$. In this case, we have 
\[f(z)=a(z)+dH(z)e^{c_1z},\]
where $d\in\mathbb{C}\backslash \{0\}$ and $H$ is transcendental such that $\rho(H)<1$ and (\ref{XMP1}) holds.

\medskip
{\bf Sub-case 1.2.} Suppose $\deg(Q)=1$. Let $Q(z)=d_1z+d_0(d_1\neq 0)$. 

We consider following sub-cases.

\medskip
{\bf Sub-case 1.2.1.} Let $c_1+d_1=0$. Then from (\ref{m1}), we have 
\bea\label{a.2} a^{(k)}-b-H_ne^{c_0+d_0}-\Big(\sideset{}{^n_{j=0}}{\sum}a_ja(z+jc)-b\Big)e^{Q}+P_k(H)e^{P}=0.\eea

Now using Lemma \ref{l3} to (\ref{a.2}), we get $P_k(H)=0$, i.e., $\sum_{j=0}^k C^k_j c_1^{k-j}H^{(j)}=0$, which shows that $H$ is a transcendental entire function and so from (\ref{se.5}), we have
\bea\label{a.2eee}|c_1^k|\leq \sideset{}{^k_{j=1}}{\sum} \big|C^k_j c_1^{k-j}\big| r^{j(\rho(H)-1+\varepsilon)}.\eea

Since $\rho(H)<1$, taking $r\to\infty$ into (\ref{a.2eee}), we get a contradiction.\par

\medskip
{\bf Sub-case 1.2.2.} Let $c_1=d_1$. Then from (\ref{m1}), we have 
\bea\label{a.3} (a^{(k)}-b)e^{-P}-H_ne^{Q}-\Big(\sideset{}{^n_{j=0}}{\sum}a_ja(z+jc)-b\Big)e^{d_0-c_0}+P_k(H)=0.\eea

Using Lemma \ref{l3} to (\ref{a.3}), we get $H_n=0$, $a^{(k)}=b$ and 
\bea\label{pp.2} \sideset{}{_{j=0}^k}{\sum} C^k_j c_1^{k-j} H^{(j)}(z)=P_k(H(z))=e^{d_0-c_0}\Big(\sideset{}{^n_{j=0}}{\sum}a_ja(z+jc)-a^{(k)}(z)\Big).\eea

Now $H_n=0$ implies that
\bea\label{pp.1}\sideset{}{^n_{j=0}}{\sum}a_je^{jcc_1}H(z+jc)=0.\eea

If $a_j=0$ for $j=0,1,\cdots, n-1$, then from (\ref{pp.1}), we deduce that $H(z+nc)=0$, which is impossible. Hence $a_j\neq 0$ for at least one $j=0,1,\cdots, n-1$. For the sake of simplicity we may assume that $a_0\neq 0$. By $\rho(H)<1$ and Lemma \ref{l8}, we obtain from (\ref{pp.1}) that $H$ is a polynomial.
Again by a simple calculation on (\ref{pp.1}), we get
\[\sideset{}{_{j=0}^n}{\sum}a_je^{jcc_1}=0.\]

If possible suppose $a=b$. Then since $a^{(k)}=b$, we have $a^{(k)}=a$, where $\rho(a)<1$. If $a$ is transcendental, then by Lemma \ref{em1}, we get a contradiction. Hence $a$ is a polynomial. Now $a^{(k)}=a$ implies that $a=0$ and so $a=b=0$. Finally from (\ref{se.3}) and (\ref{pp.2}), we get $H=0$, which is impossible. Hence $a\neq b=a^{(k)}\not\equiv 0$.

\medskip
First we suppose $a(\not\equiv 0)$ is a polynomial. Then by a simple calculation on (\ref{pp.2}), we get $\deg(H)=\deg(a)$.\par

\smallskip
Next we suppose $a$ is a transcendental entire function. 
Since $a$ is a transcendental entire function of order less than $1$, we get 
\bea\label{sm.1} P_k(H(z))/a(z)\to 0\;\;\text{and}\;\;a^{(k)}(z)/a(z)\to 0\;\text{as}\;z\to\infty.\eea

Also using Lemma \ref{em1}, we get 
\bea\label{sm.2}\sideset{}{_{j=0}^n}{\sum}a_ja(z+jc)/a(z)\to \sideset{}{_{j=0}^n}{\sum} a_j.\eea

Therefore using (\ref{sm.1}) and (\ref{sm.2}) to (\ref{pp.2}), we get $\sum_{j=0}^n a_j=0$. Since $\sum_{j=0}^n a_je^{jcc_1}=0$, we deduce that $e^{cc_1}=1$. 

\medskip
{\bf Sub-case 1.2.3.} Let $c_1\neq -d_1$ and $c_1\neq d_1$. Then from (\ref{m1}), we get 
 \bea\label{a.4} a^{(k)}-b-\Big(\sideset{}{^n_{j=0}}{\sum}a_ja(z+jc)-b\Big)e^{Q}+P_k(H)e^{P}-H_ne^{P+Q}=0.\eea

Again using Lemma \ref{l3} to (\ref{a.4}), we get $P_k(H)=0$. In this case also, we immediately get a contradiction.\par

\medskip
{\bf Case 2.} Suppose $\rho(f)\geq 2$. Since $\rho(a),\rho(b)<\max\{1,\rho(f)\}$, we have $\rho(a), \rho(b)<\rho(f)$ and so $\rho(f-a)=\rho(f)=\rho(f-b)$. By the given condition, we have $\lambda(f-a)<\rho(f)=\rho(f-a)$. This shows that $0$ and $\infty$ are the Borel exceptional values of $f-a$ and so $f-a$ is a function of regular growth, i.e., $\mu(f-a)=\rho(f-a)=\rho(f)$. Therefore $\rho(a), \rho(b)<\rho(f)=\mu(f-a)$ and so by Theorem 1.18 \cite{YY1}, we can say that both $a$ and $b$ are small functions of $f-a$. Hence both $a$ and $b$ are small functions of $f$. Now by Lemma \ref{l2.1}, we have $\deg(Q)=\rho(f)-1$. Since $\rho(f)\geq 2$, we have $\deg(Q)\geq 1$. 

If $P_k(H)-H_ne^Q\not=0$, then we see that the order of the left side of (\ref{m1}) is less than $\rho(f)$, but the order of the right side of (\ref{m1}) is equal to $\rho(f)$, which is impossible. Hence $P_k(H)-H_ne^Q=0$ ans so from (\ref{a1.a1}), we have
\bea\label{al.24} \sideset{}{^n_{j=1}}{\sum} a_je^{R_j(z)}H(z+jc)/H(z)+a_0=e^{-Q(z)}P_k(H(z))/H(z),\eea
where $R_j(z)=P(z+jc)-P(z)\;(j=1,\ldots,n)$. Since $P_k(H)-H_ne^Q=0$, from (\ref{m1}), we have 
\bea\label{xe.1} (a^{(k)}-b)e^{-Q}=\sideset{}{^n_{j=0}}{\sum}a_ja(z+jc)-b.\eea

Let 
\bea\label{al.23a} -Q(z)=d_{s-1}z^{s-1}+d_{s-2}z^{s-2}+\cdots+d_0.\eea

Then from (\ref{al.1a}), we may assume that
\bea\label{al.24a} R_j(z)=jsc_scz^{s-1}+P_{s-2,j}(z),\eea
where $\deg(P_{s-2,j})\leq s-2$. Then from (\ref{al.23a}) and (\ref{al.24a}), we get
\[R_j(z)+Q(z)=(jsc_sc-d_{s-1})z^{s-1}+\ldots,\; j=1,2,\ldots,n.\]

Now we divide following two sub-cases:\par

\smallskip
{\bf Sub-case 2.1.} Suppose there exists $j_0(1\leq j_0\leq n)$ such that $j_0sc_sc=d_{s-1}$. Therefore $\deg(R_{j_0}+Q)\leq s-2$.

Now we divide following two sub-cases:\par

\smallskip
{\bf Sub-case 2.1.1.} Suppose $n=1$. Then from (\ref{al.24}), we have
\bea\label{al.25}\label{uu.1} -a_0e^{-R_1(z)}=a_1H(z+c)/H(z)-e^{-R_1(z)-Q(z)}P_k(H(z))/H(z).\eea

\smallskip
First we suppose $a_0\neq 0$. 
Note that 
\[\mu(e^{R_1})=\rho(e^{R_1})=\deg(R_1)=s-1,\;\rho(e^{-Q-R_{1}})=\deg(-Q-R_{1})\leq s-2<s-1=\mu(e^{R_1}).\]

So $T(r,e^{-Q-R_{1}})=S(r,e^{R_1}).$ 
It is easy to prove from (\ref{e.6a}) that 
\[m(r,H(z+c)/H(z))=S(r,e^{R_1}).\] 

We know that $m(r,P_k(H)/H)=O(\log r)$ and so 
\[m(r,P_k(H)/H)=S(r,e^{R_1}).\]

On the other hand, from (\ref{uu.1}), we get 
\[m(r,e^{R_1})=m(r,e^{-R_1})=S(r,e^{R_1}).\]

Therefore $T(r,e^{R_1})=m(r,e^{R_1})=S(r,e^{R_1})$, which is a contradiction.\par 

\smallskip
Next we suppose $a_0=0$. Then from (\ref{al.25}), we get
\bea\label{uu.2} a_1H(z+c)e^{Q+R_1}=P_k(H).\eea

We consider following two sub-cases.\par

\smallskip
{\bf Sub-case 2.1.1.1.} Suppose $\rho(H)<1$. Note that $\rho(H(z+c))<1$ and $\rho(P_k(H))<1$. Then from (\ref{uu.2}), we get $\rho(e^{Q+R_1})<1$ and so $e^{Q+R_1}$ is a constant. Let $a_1e^{Q+R_1}=C_1$. So from (\ref{uu.2}), we get
\bea\label{uu.3} C_1H(z+c)=P_k(H(z)).\eea

Now from (\ref{e.2aaa}), we see that $p_i$ is a non-zero polynomial such that $\deg(p_i)=(k-i)\deg(P^{(1)})$. Then using Lemma \ref{l6a} (taking $\varepsilon=\frac{1}{2}$), we get
\bea\label{xx.0} |A_i|r^{(k-i)(s-1)}\leq2\left|p_{i}(z)\right|\leq 3|A_i|r^{(k-i)(s-1)},\eea
where $A_i$ is the leading coefficient of the polynomial $p_i$.

\smallskip
First suppose $H$ is a polynomial. 

Then from (\ref{e.2aa}), we see that $P_k(H)$ is a non-constant polynomial such that $\deg(H)<\deg(P_k(H))$. Clearly $\deg(H(z+c))<\deg(P_k(H(z)))$ and so from (\ref{uu.3}), we get a contradiction. 

\smallskip
Next suppose $H$ is a transcendental entire function. 
  
Then using (\ref{se.5}) and (\ref{xx.0}), we get
\bea\label{ww.1} &&\left|P_k(H(z))/H(z)\right|\\&=&\left|H^{(k)}(z)/H(z)+p_{k-1}(z)H^{(k-1)}(z)/H(z)+\ldots+p_0(z)\right|\nonumber\\&\geq&
|p_0(z)|-|p_1(z)|\left|H^{(1)}(z)/H(z)\right|-\ldots-\left|H^{(k)}(z)/H(z)\right|\nonumber\\&\geq&
|p_0(z)|\left\lbrace 1-|p_1(z)|/|p_0(z)|\;|z|^{\rho(H)-1+\varepsilon}-\ldots-1/|p_0(z)|\;|z|^{k(\rho(H)-1+\varepsilon)}\right\rbrace\nonumber\\&\geq&
\frac{|A_0|}{2}|z|^{k(s-1)}\left\lbrace 1-|p_1(z)|/|p_0(z)|\;|z|^{\rho(H)-1+\varepsilon}-\ldots-1/|p_0(z)|\;|z|^{k(\rho(H)-1+\varepsilon)}\right\rbrace.\nonumber\eea

As $\deg(p_i)<\deg(p_{i+1})$ and $\rho(H)<1$, we get
\[|p_1(z)|/|p_0(z)|\;|z|^{\rho(H)-1+\varepsilon}+\ldots+1/|p_0(z)|\;|z|^{k(\rho(H)-1+\varepsilon)}<\varepsilon\]
for sufficiently large values of $|z|=r$. Then from (\ref{ww.1}), we get
\bea\label{ww.2}\left|P_k(H(z))/H(z)\right|\geq |A_0|(1-\varepsilon)r^{k(s-1)}/2\eea
for sufficiently large values of $r$. Now using Lemma \ref{em1}, we conclude from (\ref{uu.3}) and (\ref{ww.2}) that
$|A_0|(1-\varepsilon)r^{k(s-1)}\leq 2|C_1|$ for sufficiently large values of $r$, which is impossible.\par

\medskip
{\bf Sub-case 2.1.1.2.} Suppose $\rho(H)\geq 1$. In this case, we have $f^{(k)}-b=e^{Q}(a_1f(z+c)-b)$ and $f(z)=a+He^{P}$, where $P$ and $Q$ are polynomials and $H$ is a transcendental entire function such that $\deg(P)\geq 2$, $\deg(Q)=\deg(P)-1$, $\rho(H)\geq 1$ and $(a^{(k)}-b)e^{-Q}=a_1a(z+c)-b$. If $a^{(k)}(z)=a_1a(z+c)$ and $a^{(k)}\not=b$, then from (\ref{xe.1}), we get $(a^{(k)}-b)(e^Q-1)=0$. Since $\deg(Q)=\rho(f)-1\geq 1$, we get a contradiction.\par

\medskip
{\bf Sub-case 2.1.2.} Suppose $n\geq 2$. 

\medskip
First we suppose $a_0\neq 0$. In this case from (\ref{al.24}), we have 
\bea\label{al.26}\left(\sum\limits_{\substack{1\leq j\leq n,j\neq j_0}}a_je^{P(z+jc)-P(z+c)}H(z+jc)/H(z)+B_{j_0}e^{P(z+j_0c)-P(z+c)}\right)e^{R_1(z)}=-a_0,\eea
where 
\bea\label{al.27} B_{j_0}(z)=a_{j_0}H(z+j_0c)/H(z)-e^{-Q(z)-R_{j_0}(z)}P_k(H(z))/H(z).\eea

Let $Q_1=e^{R_1}$. Note that 
\[Q_1(z+(j-1)c)\ldots Q_1(z+c)=e^{P(z+jc)-P(z+c)},\]
$j=2,3,\ldots,n.$ Then (\ref{al.26}) can be written as
\bea\label{al.28} U(z,Q_1)Q_1=-a_0,\eea
where 
\beas U(z,Q_1(z))&=&\sum\limits_{\substack{1\le j\leq n,j\neq j_0}}a_{j}Q_1(z+(j-1)c)Q_1(z+(j-2)c)\cdots Q_1(z+c)H(z+jc)/H(z)\\
 &&+B_{j_0}(z)Q_1(z+(j_0-1)c)Q_1(z+(j_0-2)c)\cdots Q_1(z+c)\eeas
if $j_0\geq 2$ and 
\beas U(z,Q_1(z))&=&\sideset{}{_{2\le j\leq n}}{\sum}a_jQ_1(z+(j-1)c)Q_1(z+(j-2)c)\cdots Q_1(z+c)H(z+jc)/H(z)\\&&+B_{j_0}(z),\;\text{if}\; j_0=1.\eeas
 
Now from (\ref{al.28}), we get $U(z,Q_1)\not\equiv 0$ and $\deg(U(z,Q_1))=n-1\geq 1$. Now we prove that if $a_{\lambda}$ is a coefficient of $U(z,Q_1)$, then $m(r,a_\lambda)= S(r,Q_1)$. Note that 
\[\mu(e^{R_1})=\rho(e^{R_1})=\deg(R_1)=s-1\;\text{and}\;\rho(e^{-Q-R_{j_0}})=\deg(Q+R_{j_0})\leq s-2<s-1=\mu(e^{R_1}).\]

So $T(r,e^{-Q-R_{j_0}})=S(r,Q_1)$. Also (\ref{e.6a}) gives 
\[m(r,H(z+jc)/H(z))=S(r,Q_1)\;(j=1,\ldots,n).\]

Since $m(r,P_k(H)/H)=S(r,Q_1)$, (\ref{al.27}) gives
\[m(r,B_{j_0}(z))\leq m\big(r,H(z+j_0c)/H(z)\big)+m\big(r,P_k(H)/H\big)+m\big(r,e^{-Q-R_{j_0}}\big)\leq S(r,Q_1).\]

Then in view of Remark \ref{r1} and using Lemma \ref{l6} to (\ref{al.28}), we get  $m(r,Q_1)=S(r,Q_1)$ and so $T(r,Q_1)=S(r,Q_1)$, which is a contradiction.\par

\medskip
Next we suppose $a_0=0$. For the sake of simplicity we may assume that $a_1\neq 0$. In this case, from (\ref{al.24}), we have 
\beas\label{al.26aa}\left(\sum\limits_{\substack{2\leq j\leq n,j\neq j_0}}a_je^{P(z+jc)-P(z+c)}H(z+cj)/H(z)+B_{j_0}e^{P(z+j_0c)-P(z+c)}\right)e^{R_1(z)}=-a_1,\eeas
where 
\[B_{j_0}(z)=a_{j_0}H(z+j_0c)/H(z)-e^{-Q(z)-R_{j_0}(z)}P_k(H(z))/H(z).\]

Now proceeding in the same way as done above, we get a contradiction.\par

\smallskip
{\bf Sub-case 2.2.} Suppose $jsc_sc\neq d_{s-1}$ for $1\leq j\leq n$. In this case from (\ref{al.24}), we get 
\bea\label{al.31} e^{-Q(z)}P_k(H(z))/H(z)=e^{d_{s-1}z^{s-1}}e^{\tilde{P}_{s-2}(z)}=\sideset{}{_{j=0}^n}{\sum} a_je^{R_j(z)}H(z+jc)/H(z),\eea
where 
\[\tilde{P}_{s-2}(z)=-Q(z)-d_{s-1}z^{s-1}=d_{s-2}z^{s-2}+d_{s-3}z^{s-3}+\cdots+d_0.\]

Again from (\ref{al.24a}) and (\ref{al.31}), we have
\bea\label{al.33} e^{-Q(z)}P_k(H(z))/H(z)&=&e^{d_{s-1}z^{s-1}}e^{\tilde{P}_{s-2}(z)}P_k(H(z))/H(z)\\&=&
\sideset{}{_{j=1}^n}{\sum} a_je^{jsc_scz^{s-1}} e^{P_{s-2,j}(z)}H(z+jc)/H(z)+a_0.\nonumber\eea

Now we consider following two sub-cases.\par

\smallskip
{\bf Sub-case 2.2.1.} Let $\rho(H)<1$.  If $H$ is a transcendental entire function, then from (\ref{ww.2}), we get 
\[\big|P_k(H(z))/H(z)\big|\geq \frac{|A_0|}{2}(1-\varepsilon)r^{k(s-1)},\;\;\text{i.e.,}\;\;|H(z)|/|P_k(H(z))|\leq \frac{2}{|A_0|(1-\varepsilon)r^{k(s-1)}}<B_1\]
for sufficiently large values of $r$, where $B_1>0$ is a real number.
Again if $H$ is a non-zero polynomial, then from (\ref{e.2aa}), we see that $P_k(H)$ is also a non-zero polynomial such that $\deg(H)<\deg(P_k(H))$. Consequently $\big|H(z)/P_k(H(z))\big|=O(1)$ as $|z|\to\infty$ and so, we let $\big|H(z)/P_k(H(z))\big|<B_2$ for sufficiently large values of $|z|$, where $B_2>0$.
Therefore in either case, we have 
\bea\label{wwx.2}\Big|H(z)/P_k(H(z))\Big|<B\;\;\text{for sufficiently large values of}\;|z|.\eea

Note that $ns|c_sc|>(n-1)s|cc_s|>\cdots>s|c_sc|$ and either $|d_{s-1}|\in\{js|c_sc|: j=1,2,\ldots,n\}$ or $|d_{s-1}|\not\in\{js|c_sc|: j=1,2,\ldots,n\}$.
Therefore if we compare $|d_{s-1}|$ with $ns|c_sc|$, $(n-1)s|c_sc|$, $\cdots$, $s|c_sc|$, then it is enough to compare $|d_{s-1}|$ with $ns|c_sc|$.

Now we divide following two sub-cases.\par

\medskip
{\bf Sub-case 2.2.1.1.} Let $ns|c_sc|\leq |d_{s-1}|$. Let $\arg (d_{s-1})=\theta_1$ and $\arg (c_sc)=\theta_2$. Take $\theta_0$ such that $\cos ((s-1)\theta_0+\theta_1)=1$. Then using Lemma \ref{l7}, we see that for any given $\varepsilon\;(0<\varepsilon<s-\rho(H))$, there exists a set $E\subset (1,\infty)$ of finite logarithmic measure such that for all $z=re^{i\theta_0}$ satisfying $|z|=r\not\in [0,1]\cup E$, we have
\bea\label{al.34} \exp(-r^{\rho(H)-1+\varepsilon})\leq \big|H(z+jc)/H(z)\big|\leq \exp(r^{\rho(H)-1+\varepsilon})\;(j=1,2,\ldots,n).\eea

Note that 
\bea\label{al.35} &&|\exp(d_{s-1}z^{s-1})|\\&=&|\exp(|d_{s-1}|r^{s-1}(\cos((s-1)\theta_0+\theta_1))+i\sin((s-1)\theta_0+\theta_1)))|\nonumber=\exp(|d_{s-1}|r^{s-1}).\nonumber\eea

Similarly, we can show that
\bea\label{al.36} |\exp(jsc_scz^{s-1})|=\exp(js|c_{s}c|r^{s-1}\cos((s-1)\theta_0+\theta_2)),\;j=1,2,\ldots,n.\eea

Suppose $\deg(P_{s-2,j})>0$. Set $P_{s-2,j}=c_{0,j}+c_{1,j}z+\ldots+c_{s-2,j}z^{s-2}$, $c_{k,j}=\varrho_{k,j}(\cos \alpha_{k,j}+i\sin\alpha_{k,j})$. Then
$P_{s-2,j}(z)=\sum_{k=0}^{s-2}\varrho_{k,j}r^{k}\{\cos(\alpha_{k,j}+k\theta_0)+i\sin(\alpha_{k,j}+k\theta_0)\}$ and so
\bea\label{ax.1} |\exp(P_{s-2,j}(z))|=\exp\Big(\sideset{}{^{s-2}_{k=0}}{\sum}\varrho_{k,j}r^{k}\cos(\alpha_{k,j}+k\theta_0)\Big).\eea

Note that 
\beas \Big|\sideset{}{^{s-2}_{k=0}}{\sum}\varrho_{k,j}r^{k}\cos(\alpha_{k,j}+k\theta_0)\Big|\leq \sideset{}{_{k=0}^{s-2}}{\sum}\varrho_{k,j}r^{k}=\varrho_{s-2,j}r^{s-2}\Big(1+\frac{\varrho_{s-3,j}}{\varrho_{s-2,j}}\frac{1}{r}+\ldots+\frac{\varrho_{0,j}}{\varrho_{s-2,j}}\frac{1}{r^{s-2}}\Big).\eeas

Since 
\[\frac{\varrho_{s-3,j}}{\varrho_{s-2,j}}\frac{1}{r}+\ldots+\frac{\varrho_{0,j}}{\varrho_{s-2,j}}\frac{1}{r^{s-2}}\rightarrow 0\] 
as $r\rightarrow \infty$, (\ref{ax.1}) gives 
\[\left|\exp\left(P_{s-2,j}(z)\right)\right|<\exp\left(|c_{s-2,j}|r^{s-2}(1+\varepsilon_0)\right)\]
for $r>r(\varepsilon_0)$ and $\varepsilon_0\in (0,1)$. Let $|C_0|=\max\{|c_{s-2,j}|: j=1,\ldots,n\}$ and $\varepsilon_0=\frac{1}{2}$. Then
\bea\label{al.37} |\exp(P_{s-2,j}(z))|<\exp\big(3|C_0|r^{s-2}/2\big),\;j=1,2,\ldots,n\eea
for sufficiently large values of $r$. If $P_{s-2,j}$ is a constant, then (\ref{al.37}) also holds. Similarly we can prove that 
\bea\label{al.38} |\exp(-\tilde P_{s-2}(z))|<\exp\big(3|d_{s-2}|r^{s-2}/2\big)\eea
for sufficiently large values of $r$. Clearly from (\ref{al.34}), (\ref{al.36}) and (\ref{al.37}), we get 
\bea\label{al.39} &&\Big|e^{jsc_scz^{s-1}} e^{P_{s-2,j}(z)}H(z+jc)/H(z)\Big|\\&\leq& \exp\Big(js|c_{s}c|r^{s-1}\cos((s-1)\theta_0+\theta_2))+r^{\rho(H)-1+\varepsilon}+3|C_0|r^{s-2}/2\Big)=\exp (jx+y),\nonumber\eea
where $x=s|c_{s}c|r^{s-1}\cos((s-1)\theta_0+\theta_2))$, $y=r^{\rho(H)-1+\varepsilon}+\frac{3|C_0|}{2}r^{s-2}$ and $j=1,2,\ldots,n$.
Now from (\ref{al.33}), (\ref{wwx.2}), (\ref{al.35}), (\ref{al.38}) and (\ref{al.39}), we have
\bea\label{al.40} \exp(|d_{s-1}|r^{s-1})&=&|\exp(d_{s-1}z^{s-1})|
\\&=&\Big|\frac{\exp(-Q(z))P_k(H(z))/H(z)}{\exp(\tilde{P}_{s-2}(z))P_k(H(z))/H(z)}\Big|\nonumber\\
&\leq & B\Big|\sideset{}{^n_{j=1}}{\sum} a_je^{jsc_scz^{s-1}} e^{P_{s-2,j}(z)}H(z+jc)/H(z)+a_0\Big|\;|\exp(-\tilde{P}_{s-2}(z))|\nonumber\\&\leq &
B\Big(\sideset{}{^n_{j=1}}{\sum} |a_j|\exp \left(jx+y\right)+|a_0|\Big)\exp\big(3|d_{s-2}|r^{s-2}/2\big)
\nonumber.\eea

\medskip
First we suppose $x\leq 0$, i.e., $\cos((s-1)\theta_0+\theta_2))\leq 0$. Let $\cos((s-1)\theta_0+\theta_2))=-x_1$, where $0\leq x_1\leq 1$. If $x_1=0$, then from (\ref{al.40}), one can easily conclude that
\bea\label{ax.2} \exp(|d_{s-1}|r^{s-1})<C_1\exp\big(r^{\rho(H)-1+\varepsilon}+3(|C_0|+|d_{s-2}|)r^{s-2}/2\big),\eea
where $C_1>0$. Since $\rho(H)-1+\varepsilon<s-1$, (\ref{ax.2}) leads to a contradiction. Next we suppose $0<x_1\leq 1$.
In this case, since $\rho(H)-1+\varepsilon<s-1$, we see that
\[\sideset{}{^n_{j=1}}{\sum} |a_j|\exp \big(-js|c_{s}c|x_1r^{s-1}+r^{\rho(H)-1+\varepsilon}+3|C_0|r^{s-2}/2\big)\rightarrow 0\] 
as $r\rightarrow \infty$. 
Now for $\varepsilon_1$ in $(0,1)$ and $r>r(\varepsilon_1)$, we get from (\ref{al.40})
\beas \exp(|d_{s-1}|r^{s-1})<(1+\varepsilon_1)\exp\big(3|d_{s-2}|r^{s-2}/2\big),\eeas
which leads to a contradiction.

\medskip
Next we suppose $x>0$, i.e., $\cos((s-1)\theta_0+\theta_2))>0$. In this case, we have $\exp (jx+y)<\exp (nx+y)$ for $j=0,\ldots,n$ and so from (\ref{al.40}), we get
\bea\label{ax.4} &&\exp(|d_{s-1}|r^{s-1})\\&\leq &
B\big(\sideset{}{^n_{j=1}}{\sum} |a_j|\exp \left(jx+y\right)+|a_0|\big)\exp\big(3|d_{s-2}|r^{s-2}/2\big)
\nonumber\\&\leq&
B\big(\sideset{}{^n_{j=0}}{\sum} |a_j|+1\big)\exp (nx+y)\exp\big(3|d_{s-2}|r^{s-2}/2\big)
\nonumber\\
&=&B_3\exp\Big(ns|c_{s}c|r^{s-1}\cos((s-1)\theta_0+\theta_2))+r^{\rho(H)-1+\varepsilon}+3(|C_0|+|d_{s-2}|)r^{s-2}/2\Big)\nonumber,\eea
where $B_3=B\big(\sum_{j=0}^n|a_j|+1\big)$. Since $\rho(H)+\varepsilon<s$ and $B_3=o(r^{s-1})$, from (\ref{ax.4}), we get
\bea\label{al.41} \exp(|d_{s-1}|r^{s-1})\leq \exp(ns|c_sc|\cos((s-1)\theta_0+\theta_2)r^{s-1}+o(r^{s-1})).\eea

By assumption, we have $d_{s-1}\neq nsc_sc$ and $ns|c_sc|\leq |d_{s-1}|$. If $ ns|c_sc|=|d_{s-1}|$, then $\cos((s-1)\theta _0+\theta _2)\neq 1$ and so $\cos((s-1)\theta_0+\theta_2)<1$. Therefore $ns|c_sc|\cos((s-1)\theta_0+\theta_2)<ns|c_sc|=|d_{s-1}|$. If $ns|c_sc|<|d_{s-1}|$, then $ns|c_sc|\cos((s-1)\theta_0+\theta_2)\leq ns|c_sc|<|d_{s-1}|$. So in either case, we have 
\[ns|c_sc|\cos((s-1)\theta_0+\theta_2)<|d_{s-1}|.\]

Now there exists $\varepsilon_2>0$ such that $ns|c_sc|\cos((s-1)\theta_0+\theta_2)+2\varepsilon_2<|d_{s-1}|$ and so from (\ref{al.41}), we have 
$\exp\left(|d_{s-1}|r^{s-1}\right)\leq \exp\left(\left(|d_{s-1}|-2\varepsilon_2\right)r^{s-1}\right)$,
which is a contradiction.\par

\medskip
{\bf Sub-case 2.2.1.2.} Let $ns|c_sc|>|d_{s-1}|$. In this case, (\ref{al.33}) can be rewritten as 
\bea\label{xa.1} &&e^{nsc_scz^{s-1}} e^{P_{s-2,n}(z)}\\&=&\sideset{}{^{n-1}_{j=0}}{\sum}\frac{a_j}{a_n}\frac{H(z+jc)}{H(z+nc)}e^{jsc_scz^{s-1}} e^{P_{s-2,j}(z)}-
\frac{1}{a_n}\frac{P_k(H(z))}{H(z+nc)}e^{d_{s-1}z^{s-1}}e^{\tilde{P}_{s-2}(z)}\nonumber.\eea

Let $\arg (c_sc)=\theta_1$ and $\arg (d_{s-1})=\theta_2$. Take $\theta_0$ such that $\cos ((s-1)\theta_0+\theta_1)=1$. Then using Lemma \ref{l7}, we see that for any given $\varepsilon\;(0<\varepsilon<s-\rho(H))$, there exists a set $E\subset (1,\infty)$ of finite logarithmic measure such that for all $z=re^{i\theta_0}$ satisfying $|z|=r\not\in [0,1]\cup E$, we have
\bea\label{xa.2} \exp(-r^{\rho(H)-1+\varepsilon})\leq \big|H(z+jc)/H(z+nc)\big|\leq \exp(r^{\rho(H)-1+\varepsilon})\;(j=0,1,\ldots,n-1).\eea

In this case, we easily get
\bea\label{xa.3} |\exp(jsc_scz^{s-1})|=\exp(js|c_{s}c|r^{s-1})\eea
$j=0,1,\ldots,n$ and 
\bea\label{xa.4} |\exp(d_{s-1}z^{s-1})|=|\exp(|d_{s-1}|r^{s-1}\cos((s-1)\theta_0+\theta_2))|\leq |\exp(|d_{s-1}|r^{s-1})|.\eea

Also we easily obtain
\bea\label{xa.5} |\exp (\pm P_{s-2,j}(z))|<\exp\big(3|C_0|r^{s-2}/2\big)\;\text{and}\;|\exp (\tilde{P}_{s-2}(z))|<\exp\big(3|d_{s-2}|r^{s-2}/2\big)\eea
for sufficiently large values of $r$. Then from (\ref{xa.2}), (\ref{xa.3}) and (\ref{xa.5}), we get
\bea\label{xa.6} &&\Big|e^{jsc_scz^{s-1}} e^{P_{s-2,j}(z)}H(z+jc)/H(z+nc)\Big|\\&\leq& \exp\big(js|c_{s}c|r^{s-1}+r^{\rho(H)-1+\varepsilon}+3|C_0|r^{s-2}/2\big)\nonumber\\&\leq & \exp\big((n-1)s|c_{s}c|r^{s-1}+r^{\rho(H)-1+\varepsilon}+3|C_0|r^{s-2}/2\big),\nonumber\eea
$j=1,\ldots, n-1$. Let $\varepsilon_3>0$ be such that $\rho(H)+\varepsilon_3<s$. Using (\ref{se.5}) and (\ref{xx.0}), we get
\beas \big|P_k(H(z))/H(z)\big|&\leq& \sideset{}{^k_{j=0}}{\sum} |p_j(z)|\big|H^{(j)}(z)/H(z)\big|\\&\leq& \sideset{}{^k_{j=0}}{\sum} 3|A_j|r^{k(s-1)+j(\rho(H)-s+\varepsilon_3)/2}\\&<&C_2r^{k(s-1)},\eeas
where $C_2=\sum_{j=0}^k \frac{3|A_j|}{2}$. For sufficiently large $r$, we get $C_2r^{k(s-1)}<\exp(r^{s-\frac{3}{2}})$ and so
\bea\label{xa.6a} \big|P_k(H(z))/H(z)\big|<\exp(r^{s-\frac{3}{2}}).\eea

Then from (\ref{xa.2}), (\ref{xa.4}), (\ref{xa.5}) and (\ref{xa.6a}), we get
\bea\label{xa.7} &&\Big|e^{d_{s-1}z^{s-1}}e^{\tilde{P}_{s-2}(z)}P_k(H(z)/H(z+nc)\Big|\\&\leq& \exp\Big(|d_{s-1}|r^{s-1}+r^{\rho(H)+s-\frac{5}{2}+\varepsilon}+3|d_{s-2}|r^{s-2}/2\Big).\nonumber
\eea

Finally from (\ref{xa.1}), (\ref{xa.3}), (\ref{xa.5})-(\ref{xa.7}), we obtain
\bea\label{xa.8} &&\exp(ns|c_{s}c|r^{s-1})\\&=&\Big|\sideset{}{^{n-1}_{j=0}}{\sum}\frac{a_j}{a_n}\frac{H(z+jc)}{H(z+nc)}e^{jsc_scz^{s-1}} e^{P_{s-2,j}(z)}-
\frac{1}{a_n}\frac{P_k(H(z))}{H(z+nc)}e^{d_{s-1}z^{s-1}}e^{\tilde{P}_{s-2}(z)}\Big||e^{P_{s-2,n}(z)}|\nonumber\\&\leq&
B_4\exp\Big(x_2r^{s-1}+r^{\rho(H)+s-\frac{5}{2}+\varepsilon}+(6|C_0|+3|d_{s-2}|)r^{s-2}/2\Big),\nonumber
\eea
where $B_4=\sum_{j=0}^n |a_j|+1$ and $x_2=\max\{(n-1)s|c_{s}c|, |d_{s-1}|\}$. Since $|d_{s-1}|<ns|c_{s}c|$, $x_2<ns|c_{s}c|$. Also $\rho(H)+s-\frac{5}{2}+\varepsilon<s-1$. Consequently from (\ref{xa.8}), we get a contradiction.\par

\medskip
{\bf Sub-case 2.2.2.} Let $\rho(H)\geq 1$. Note that $P_k(H)/H=\left(He^P\right)^{(k)}/He^P$. Let $g=He^P$. 

\medskip
First suppose $|g^{(k)}(z)|$ is bounded on the ray $\arg z=\theta_0$. By considering the formula
\bea\label{xx.1} g^{(k-1)}(z)=g^{(k-1)}(0)+\int_{0}^zg^{(k)}(w)dw,\eea
we obtain that $|g^{(k-1)}(z)|\leq |g^{(k-1)}(0)|+M|z|$ for all $z$ satisfying $\arg z=\theta_0$, where $M=\max\{|g^{(k)}(re^{i\theta_0})|\}>0$ is some constant.
Again by using the formula
\bea\label{xx.3} g^{(j-1)}(z)=g^{(j-1)}(0)+\int_{0}^zg^{(j)}(w)dw\;\;(1\leq j\leq k),\eea
we easily obtain
\bea\label{xx.4} |g^{(k-2)}(z)|\leq (M+o(1))|z|^2\eea
as $r\to +\infty$ for all $z$ satisfying $\arg z=\theta_0$. Similarly, by using (\ref{xx.4}) and the formula (\ref{xx.3}), for $k$ times, we easily obtain
\bea\label{xx.5} |g(z)|\leq (M+o(1))|z|^k\eea
as $r\to +\infty$ for all $z$ satisfying $\arg z =\theta_0$. Since $|g^{(k)}(z)|$ is bounded on the ray $\arg z=\theta_0$, there exists an infinite sequence of points $z_m=r_me^{i\theta_0}$, where $r_m\to\infty$ such that $|g^{(k)}(z_m)|>0$. Therefore there exists a positive constant $m$ such that 
$0<m\leq |g^{(k)}(z_m)|<M$ and so from (\ref{xx.5}), we get
\beas\label{xx.7} \big|H(z_m)/P_k(H(z_m))\big|=\big|g(z_m)/g^{(k)}(z_m)\big|\leq \Big(\frac{M}{m}+o(1)\Big)|z_m|^k=\Big(\frac{M}{m}+o(1)\Big)r_m^k\;\text{as}\;r_m\to\infty.\eeas
 
\medskip
Next suppose $|g^{(k)}(z)|$ is unbounded on the ray $\arg z=\theta_0$. Then by Lemma \ref{l2a}, there exists an infinite sequence of points $z_m=r_me^{i\theta_0}$, where $r_m\to\infty$ such that $g^{(k)}(z_m)\to\infty$ and $|g(z_m)|\leq (1+o(1))|z_m|^k|g^{(k)}(z_m)|$ as $r_m\to \infty$ and so 
\beas\label{yy.1} \big|H(z_m)/P_k(H(z_m))\big|=\big|g(z_m)/g^{(k)}(z_m)\big|\leq (1+o(1))|z_m|^k=(1+o(1))r_m^k\;\text{as}\;r_m\to\infty.\eeas

Thus in either case, we have 
\bea\label{xx.7} \big|H(z_m)/P_k(H(z_m))\big|=\big|g(z_m)/g^{(k)}(z_m)\big|\leq \Big(\frac{M}{m}+o(1)\Big)|z_m|^k=\Big(\frac{M}{m}+o(1)\Big)r_m^k\eea
as $r_m\to\infty$. 
Note that $r_m^k<\exp(r_m^{\rho(H)-1+\varepsilon})$ for sufficiently large $m$. Now from (\ref{xx.7}), we get
\bea\label{xx.8} \big|H(z_m)/P_k(H(z_m))\big|<\exp(r_m^{\rho(H)-1+\varepsilon})\eea 
for sufficiently large values of $m$.

If $ns|c_sc|\leq |d_{s-1}|$, then from (\ref{al.33}), (\ref{al.35}), (\ref{al.38}), (\ref{al.39}) and (\ref{xx.8}), we have
\beas &&\exp(|d_{s-1}|r_m^{s-1})\\&=&
\left|\frac{\exp(-Q(z))P_k(H(z_m))/H(z_m)}{\exp(\tilde{P}_{s-2}(z_m))P_k(H(z_m))/H(z_m)}\right|\nonumber\\
&\leq & \left|\sideset{}{^n_{j=1}}{\sum}a_j\frac{H(z_m+jc)}{H(z_m)}e^{jsc_scz_m^{s-1}} e^{P_{s-2,j}(z_m)}+a_0\right|\;\left|\frac{H(z_m)}{P_k(H(z_m))}\right|\;\left|\exp(-\tilde{P}_{s-2}(z_m))\right|\nonumber\\&\leq &
\big(\sideset{}{^n_{j=1}}{\sum} |a_j|\exp \left(jx+y\right)+|a_0|\big)\exp\Big(r_m^{\rho(H)-1+\varepsilon}+3|d_{s-2}|r_m^{s-2}/2\Big)
\nonumber.\eeas

Now proceeding in the same way as done in Sub-case 2.2.1.1, we get a contradiction. Hence $ns|c_sc|>|d_{s-1}|$. In this case also we get a contradiction from Sub-case 2.2.1.2. 

This completes the proof.
\end{proof}

\section{Open questions}

As we know $\lambda(f)\leq \rho(f)$ and since in Theorem \ref{t2.1}, we have dealt with the case $\lambda(f)< \rho(f)$, it will be interesting to inspect the existence of the solutions of the Eq. (\ref{xx}) for the case $\lambda(f)=\rho(f)$.  In the following example, we point out the fact that when $\lambda(f)=\rho(f)$, the solution of Eq. (\ref{xx}) may exist.
\begin{exm} Let $f(z)=\cos z$, $a_0,a_1,\ldots,a_n\in\mathbb{C}$ such that $a_0+a_2+\ldots+a_{n-2}+a_n=1$, $a_1+a_3+\ldots+a_{n-3}+a_{n-1}=0$, where $n$ is even, $c=\pi$ and $k=4$. Then $\lambda(f)=\rho(f)=1$ and $L_c(f(z))=\sin z$. Clearly $f(z)$ satisfies the Eq. $f^{(k)}(z)=L_c(f(z))$.
\end{exm}

Again since in Theorem \ref{t1}, we considered the case $\lambda(f-a)<\rho(f)$, so the conclusions of Theorem \ref{t1} under the case $\lambda(f-a)<\rho(f)$ is an enigma. So, here we place the following questions:
\begin{ques} Can one find the precise form of the solutions of Eq. (\ref{xx}) under the hypothesis $\lambda(f)=\rho(f)$?
\end{ques}

\begin{ques} Does Theorem \ref{t1} hold under the hypothesis $\lambda(f-a)=\rho(f)$?
\end{ques}

\section{Statements and declarations}
\vspace{1.3mm}

\noindent \textbf {Conflict of interest:} The authors declare that there are no conflicts of interest regarding the publication of this paper.\vspace{1.5mm}

\noindent{\bf Funding:} There is no funding received from any organizations for this research work.
\vspace{1.5mm}

\noindent \textbf {Data availability statement:}  Data sharing is not applicable to this article as no database were generated or analyzed during the current study.

\end{document}